\documentclass[10pt]{amsart}
\usepackage[T1]{fontenc}
\usepackage{lmodern}
\usepackage{amsmath,amssymb,amsthm,mathtools}
\usepackage[margin=3.6cm]{geometry}
\usepackage{microtype}
\usepackage[hidelinks]{hyperref}
\usepackage{bookmark}
\hypersetup{pdftitle={Nested algebraic limit cycles on distinct irreducible curves in quadratic systems},
pdfauthor={Chen Cao and Jaume Llibre},
pdfsubject={Exclusion of nested limit cycles on distinct irreducible algebraic carriers for planar quadratic systems},
pdfkeywords={quadratic differential systems, algebraic limit cycles, invariant algebraic curves, cofactors, periodic averages}}
\newcommand{\R}{\mathbb R}
\newcommand{\C}{\mathbb C}
\newcommand{\Z}{\mathbb Z}
\DeclareMathOperator{\divergence}{div}
\newtheorem{maintheorem}{Theorem}

\newtheorem{theorem}{Theorem}[section]
\newtheorem{lemma}[theorem]{Lemma}
\newtheorem{proposition}[theorem]{Proposition}
\newtheorem{corollary}[theorem]{Corollary}
\theoremstyle{definition}

\theoremstyle{remark}

\numberwithin{equation}{section}
\begin{document}
\title[Nested limit cycles on distinct algebraic curves]
{Nested algebraic limit cycles on distinct irreducible curves in quadratic systems}
\author[C. Cao]
{Chen Cao${}^{*}$}
\thanks{${}^{*}$Corresponding author.}

\address{\textsc{Chen Cao}\newline
School of Mathematical Sciences, CMA-Shanghai, Shanghai Jiao Tong University,
Shanghai, 200240, P. R. China}
\email{aceaccce@sjtu.edu.cn}

\author[J. Llibre]
{Jaume Llibre}

\address{\textsc{Jaume Llibre}\newline
Departament de Matem\`{a}tiques, Universitat Aut\`{o}noma de Barcelona,
08193 Bellaterra, Barcelona, Catalonia, Spain \newline
Reial Acad\`emia de Ci\`encies i Arts de Barcelona, La Rambla 115,
08002 Barcelona, Catalonia, Spain}
\email{jaume.llibre@uab.cat}
\date{}
\begin{abstract}
We prove that a real planar polynomial differential system of degree
at most two cannot have two nested limit cycles contained in distinct
irreducible invariant algebraic curves. There is no bound on the degrees
of the algebraic carriers, and neither hyperbolicity nor adjacency of
the cycles is assumed. Cofactor identities and a one-sided return-map
argument imply that the two cofactors and the divergence are linearly
independent as affine functions. They also determine the time averages
and characteristic exponents per unit time of the periodic orbits on
each carrier.
For a genuinely quadratic field, we prove that a periodic orbit whose
time average is an equilibrium with nonreal eigenvalues implies the
existence of a real invariant affine line. This excludes parallel
cofactor gradients and carriers through the common enclosed equilibrium.
For independent gradients, local invariant branches and the line joining
the two time averages restrict the finite equilibria and the common
directions of the highest homogeneous components. B\'ezout's theorem and
the index on a large circle then give a contradiction. The proof allows
real invariant affine lines and singular algebraic carriers. It does not
determine the number of limit cycles on a single irreducible carrier.
\end{abstract}
\maketitle

\section{Introduction and statement of the main result}
\label{sec:introduction}

Consider a real polynomial differential system
\begin{equation}\label{eq:system}
 \dot x=P(x,y),\qquad \dot y=Q(x,y),
\end{equation}
where $P,Q\in\R[x,y]$ and $\max\{\deg P,\deg Q\}\le2$.
We write $\mathcal{X}=P\partial_x+Q\partial_y$ for the associated
vector field and set $D=\divergence\mathcal{X}=P_x+Q_y$.
The notation $D\mathcal{X}(p)$ denotes its Jacobian matrix at $p$.
Throughout this paper, periodic orbits are understood to be nonconstant,
so equilibria are excluded. An isolated periodic orbit is called a
\emph{limit cycle}. It is \emph{algebraic} if it is contained in the
zero set of a nonzero real polynomial.

For a nonconstant polynomial $F\in\C[x,y]$, the algebraic curve
$F=0$ is invariant if
\begin{equation}\label{eq:invariance}
 \mathcal{X}F=K_FF
\end{equation}
for a polynomial $K_F\in\C[x,y]$, called its \emph{cofactor}.
Since system~\eqref{eq:system} has degree at most two, $K_F$ has
degree at most one.

Throughout this paper, the \emph{carrier} of an algebraic periodic orbit means its
minimal real irreducible, reduced invariant algebraic carrier. More
precisely, it is the zero set of an irreducible $F\in\R[x,y]$ vanishing
on the orbit, with $F$ taken only to the first power. Such an $F$ is
unique up to multiplication by a nonzero real constant. Indeed, one
irreducible factor of any polynomial vanishing on an analytic orbit
vanishes identically on it. That factor is absolutely irreducible, as
shown in Lemma~\ref{lem:regular}. Since $\mathcal{X}F$ also vanishes on the orbit,
the finiteness of the intersection of two coprime complex plane curves
implies $F\mid \mathcal{X}F$. The same argument proves uniqueness.
In particular, the carrier is attached to the chosen orbit and contains
no additional algebraic component.

An \emph{inverse integrating factor} is a $C^1$ function $V\not\equiv0$
defined on an open set and satisfying $\mathcal{X}V=DV$.

Two disjoint periodic orbits are nested if the closed Jordan disk
bounded by one is contained in the open disk bounded by the other.
Our main result is the following.

\begin{maintheorem}\label{thm:main}
System~\eqref{eq:system} cannot have two nested limit cycles
contained in distinct irreducible
invariant algebraic curves.
\end{maintheorem}

There is no bound on the degrees of the carriers in Theorem~\ref{thm:main}.
Hyperbolicity is not assumed, and the two cycles need not be adjacent.
The system may have real
invariant affine lines, and the carriers may have singular points,
isolated real points and other local branches inside the periodic disks.
No condition of global smoothness is imposed. The theorem concerns two
different minimal carriers; it does not establish uniqueness of algebraic
limit cycles on a single irreducible curve.

Chavarriga, Garc\'ia and Sorolla~\cite{ChavarrigaGarciaSorolla2006}
studied the non-nested configuration for two distinct irreducible
carriers. A precise statement and proof of the non-nested exclusion
appear in Sorolla's thesis~\cite[Theorem~5.2 and Hypothesis~A]{Sorolla2004Thesis}.
Theorem~\ref{thm:main} combined with the
non-nested exclusion in \cite[Theorem~5.2]{Sorolla2004Thesis} gives the following
statement about carriers only.

\begin{corollary}\label{cor:carriers}
For a real planar polynomial differential system of degree at most two,
at most one minimal real irreducible invariant algebraic carrier can
contain a limit cycle.
\end{corollary}
\begin{proof}
If two different carriers each contained a limit cycle, choose one on
each. Distinct periodic orbits are disjoint Jordan curves. Their disks
are either nested or disjoint. The former case is excluded by
Theorem~\ref{thm:main}, and the latter by the cited non-nested theorem,
whose hypotheses require two different real irreducible carriers.
\end{proof}

The distinction between the two configurations is essential to the
proof. For non-nested cycles, the two periodic disks contain distinct
foci $p_1,p_2$. If a carrier $F_i=0$, with cofactor $K_i$, avoids both
foci, then $K_i(p_1)=K_i(p_2)=0$. The projective two-point criterion
then forces all limit cycles onto that carrier, contradicting the
limit cycle on the other carrier. Since two distinct real irreducible carriers
cannot share the conjugate local branches at a focus, each carrier
must pass through exactly one of the two foci. After relabeling the
carriers, their cofactor values satisfy
\[
 K_1(p_1)=D(p_1),\qquad K_2(p_2)=D(p_2),\qquad
 K_1(p_2)=K_2(p_1)=0;
\]
see \cite[Lemmas~5.4--5.5]{Sorolla2004Thesis}. Thus the affine function
$D-K_1-K_2$ vanishes at two distinct foci. In the real saddle case,
the local branch restrictions give a third noncollinear zero, forcing
$D=K_1+K_2$. The other finite singular configurations are either
excluded or lead to the same identity; see
\cite[Propositions~5.7--5.9]{Sorolla2004Thesis}.
The identity yields the polynomial inverse integrating factor $F_1F_2$ used
in the non-nested proof.

For nested cycles, the two disks have only one common interior
equilibrium $e$. A carrier avoiding $e$ gives just $K_i(e)=0$, so the
first use of the two-point criterion above is no longer available.
Moreover, Proposition~\ref{prop:avoid-focus} below shows that a
hypothetical nested pair would have both carriers avoiding $e$, with
$K_1(e)=K_2(e)=0$ and $D(e)\ne0$. Hence $D-K_1-K_2$ is nonzero at $e$,
and the identity used to finish the non-nested argument cannot hold.
The common equilibrium therefore requires additional information
from the periodic orbits themselves.

We obtain this information from their actual time averages $m_1,m_2$.
The cofactor integrals give
\[
 K_j(m_i)=0,\qquad (D-K_i)(m_i)=0\qquad(j\ne i).
\]
We prove that these averages are distinct, and that the nonzero affine
function $S=D-K_1-K_2$ vanishes on the line joining them.
The key additional step is Lemma~\ref{lem:meanline}: for a genuinely
quadratic field, a periodic average equal to an equilibrium with
nonreal eigenvalues forces a real invariant affine line. Its proof
uses a periodic logarithmic integral as well as the covariance of the
actual orbit. Together with Zhou and Wang's counting theorem, it
excludes carriers through $e$, parallel cofactor gradients, and periodic
orbits outside the two carriers. An annulus argument then gives opposite
stabilities without assuming adjacency. The resulting signs and the
local cofactor values place any possible real saddle in the segment
joining $m_1$ and $m_2$, inside the outer periodic disk, where no such
saddle can occur. These connections between the common equilibrium,
periodic averages and local branches supply the restrictions needed
for the final B\'ezout and index argument.

The proof uses established tools from the theory of invariant algebraic
curves.
The relation between cofactors, characteristic exponents and periodic
averages, including the linear systems determining stability on two
invariant curves, is treated by Gasull and Giacomini
\cite[Section~5, especially Corollary~5.3]{GG2023}. The relation between
return-map multiplicity and inverse integrating factors is studied in
\cite{GGG2010}. We supply the one-sided argument needed for arbitrary real
powers before invoking the polynomial inverse integrating factor theorem
of Coll, Ferragut and Llibre \cite[Theorem~3]{CollFerragutLlibre2010IIF}.
The convexity and equilibrium properties of quadratic periodic disks
are classical; see \cite[p.~115]{Bamon1986} and
\cite[proof of Corollary~1.5]{GG2023}.
The local branch results of Chavarriga, Giacomini and Grau
\cite{ChavarrigaGiacominiGrau2005Local} constrain the possible cofactors
and intersections of invariant curves. Related arguments at finite
equilibria and at infinity occur in~\cite[Chapter~5]{Sorolla2004Thesis}.
We prove the versions needed here, including their spectral hypotheses.
We also give a direct proof of the two special cases of the projective
two-point criterion \cite[Theorem~2.4]{ChavarrigaGarcia2006Cubic} needed below.

\paragraph{Organization of the proof.}
Section~\ref{sec:preliminaries} collects the facts about minimal
carriers, characteristic exponents and one-sided return maps. They give
affine independence of the cofactors and divergence, and determine the
time average on each carrier. Section~\ref{sec:periodic-disks} proves
the periodic-disk properties used below. It also shows that an actual
periodic average at the enclosed equilibrium gives a real invariant line
and an explicit normal form covered by Zhou and Wang's uniqueness
theorem~\cite[Theorem~1]{ZhouWang2003Line}.

Section~\ref{sec:parallel-cofactors} excludes the parallel-gradient
case. For independent gradients, Section~\ref{ind:sec:means} uses the
cofactors as affine coordinates and proves that the two characteristic
exponents have opposite signs. The cofactor zero lines and local branch
results in Section~\ref{ind:sec:local} then restrict the finite
equilibria in Section~\ref{ind:sec:finite}. Finally,
Section~\ref{ind:sec:infinity} excludes common directions of the highest
homogeneous components and gives the B\'ezout and index contradiction.

\section{Cofactors, isolation and forced time averages}
\label{sec:preliminaries}

Suppose that there are two distinct carriers $F_1,F_2$ containing
limit cycles $\Gamma_1,\Gamma_2$. We set $K_i=K_{F_i}$, so that
$\mathcal{X}F_i=K_iF_i$ for $i=1,2$.
The results of this section do not require nesting. Let $T_i$ be the
period and put
\[
 \Lambda_i=\int_0^{T_i}D(\Gamma_i(t))\,dt,\qquad
 m_i=\frac1{T_i}\int_0^{T_i}\Gamma_i(t)\,dt.
\]
Thus $\Lambda_i$ is the characteristic exponent of $\Gamma_i$, and
its transverse multiplier is $e^{\Lambda_i}$.
We also write $\omega_{\mathcal{X}}=-Q\,dx+P\,dy$ for the polynomial
one-form associated with $\mathcal{X}$.

We may suppose that $\mathcal{X}$ has degree exactly two and that $P,Q$ are
coprime over $\C$. An affine linear planar field has no isolated
periodic orbit: any periodic orbit belongs to a center annulus.
If the real components $P,Q$ have a nonconstant common complex factor,
conjugation gives a nonconstant real common divisor $h$.
Writing $\mathcal{X}=h\mathcal{Y}$, the field $\mathcal{Y}$ has degree at most one. A nonconstant
periodic orbit avoids $h=0$, as does a sufficiently small collar
around it. On that collar $\mathcal{X}$ and $\mathcal{Y}$ have the same orbits, possibly with
reversed time, so isolation is again impossible. This is only a local
orbit equivalence away from $h=0$. Also, if $D\equiv0$, the form
$\omega_{\mathcal{X}}$ is closed and has a polynomial primitive whose
differential is nonzero on a periodic orbit. The resulting first integral
excludes isolation. These observations also remove the degree,
common-factor and Hamiltonian cases separated in
\cite[Remark~2]{CollFerragutLlibre2010IIF}.

\begin{lemma}[Minimal carriers and characteristic exponents]
\label{lem:regular}
Each carrier $F_i$ is absolutely irreducible and satisfies
$dF_i\ne0$ on $\Gamma_i$. For $j\ne i$, $F_j$ is nowhere zero in
a collar of $\Gamma_i$. Moreover,
\begin{equation}\label{eq:integrals}
 \int_0^{T_i}K_i(\Gamma_i(t))\,dt=\Lambda_i,
 \qquad \int_0^{T_i}K_j(\Gamma_i(t))\,dt=0\quad(j\ne i).
\end{equation}
These conclusions hold for every periodic orbit on either carrier,
whether or not isolation is assumed for that orbit.
\end{lemma}
\begin{proof}
If a real irreducible polynomial is not absolutely irreducible, it is
a constant multiple of a product of two distinct conjugate complex
irreducible polynomials. Its real zero set lies in their finite
intersection, and cannot contain a periodic orbit.
An irreducible complex plane curve has finitely many singular points.
For the local flow $\phi_t$, invariance gives
\[
 F_i(\phi_t(p))=F_i(p)
 \exp\left(\int_0^t K_i(\phi_s(p))\,ds\right).
\]
Differentiating on $F_i(p)=0$ shows that the singular locus is invariant
under local flow diffeomorphisms. A finite invariant set consists of
equilibria. Likewise, the finite intersection of the two distinct
absolutely irreducible invariant curves consists of equilibria.
Neither set meets a nonconstant orbit, and compactness gives the collar.

On an analytic transversal, $s=F_i$ is a valid coordinate.
Its return multiplier equals $\exp\int K_i\,dt$. The variational
equation has tangential multiplier one, so Liouville's formula gives
the transverse multiplier as $\exp\int D\,dt$. The exponents are real,
hence equal. For the other carrier one may integrate
$\mathcal{X}\log|F_j|=K_j$ on the orbit, obtaining zero. No logarithmic derivative
of the vanishing function $F_i$ is used on its own orbit.
\end{proof}

\begin{lemma}[Dependent cofactors and one-sided transverse orders]
\label{lem:weights}
The cofactors $K_1,K_2$ are linearly independent over $\R$.
If $D=aK_1+bK_2$ with $a,b\in\R$, then isolation of the two
chosen cycles forces $a,b$ to be positive integers.
\end{lemma}
\begin{proof}
If $K_i=0$, then $F_i$ is a regular analytic first integral near
$\Gamma_i$, contrary to isolation. Otherwise a dependence
$rK_1+sK_2=0$ has $rs\ne0$. Near $\Gamma_1$ the function
$F_1|F_2|^{s/r}$ is a single-valued real analytic first integral,
regular on the cycle. Its transverse coordinate is preserved by the
return map, which is therefore the identity. This also contradicts
isolation.

Now assume $D=aK_1+bK_2$. On one punctured side of a collar around
$\Gamma_1$ set $\mu=|F_1|^{-a}|F_2|^{-b}$. Then
$\divergence(\mu \mathcal{X})=0$, so the one-form
$\omega=\mu\omega_{\mathcal{X}}$ is closed and annihilates $\mathcal{X}$.
It is preserved by the flow, by Cartan's formula; the variable return
time adds no term since $\omega(\mathcal{X})=0$. On the side $s=F_1>0$ of
a transversal write
$\omega=s^{-a}u(s)\,ds$, where $u$ is analytic across zero and
$u(0)\ne0$. Here $u$ includes the transverse flux factor.
For the analytic return map $\Pi$,
\begin{equation}\label{eq:return}
 \Pi(s)^{-a}u(\Pi(s))\Pi'(s)=s^{-a}u(s),\qquad s>0.
\end{equation}
If $q=\Pi'(0)\ne1$, then $q>0$ and leading terms give
$q^{1-a}=1$, hence $a=1$.
If $q=1$, isolation and analyticity imply
\[
 \Pi(s)=s+cs^m+O(s^{m+1}),\qquad c\ne0,\quad m\ge2.
\]
The alternative of infinite order would make $\Pi$ the identity.
Since $u(\Pi(s))/u(s)=1+O(s^m)$, comparison of order $m-1$ in
\eqref{eq:return} gives $(m-a)c=0$, so $a=m$.
The same proof at $\Gamma_2$ gives the assertion for $b$.
It covers zero, negative and nonintegral trial weights, as well as
semistable isolated cycles, without asserting that the real powers
extend analytically across the cycle.
\end{proof}

\begin{proposition}[Affine independence and time averages]
\label{prop:means}
The three affine functions $K_1,K_2,D$ are linearly independent in
$\R[x,y]_{\le1}$. Every periodic orbit on $F_i=0$ is hyperbolic
and has the same time average $m_i$, uniquely determined by
\begin{equation}\label{eq:means-general}
 K_j(m_i)=0,\qquad (D-K_i)(m_i)=0,\qquad j\ne i.
\end{equation}
For any such orbit with period $T$ and characteristic exponent $\Lambda$,
\begin{equation}\label{eq:rate-general}
 D(m_i)=K_i(m_i)=\Lambda/T\ne0.
\end{equation}
Both cofactors are nonconstant, and $m_1\ne m_2$.
\end{proposition}
\begin{proof}
By Lemma~\ref{lem:weights}, any dependence of the three functions
would express $D=aK_1+bK_2$ with $a,b$ positive integers.
Then $V=F_1^aF_2^b$ is a nonzero polynomial satisfying $\mathcal{X}V=DV$.
Theorem~3 of \cite{CollFerragutLlibre2010IIF} states that a quadratic system with a polynomial
inverse integrating factor has no limit cycles. It has no degree bound
on $V$ and applies after the preliminary reductions above.
Thus the three affine functions are independent.

Affine functions commute with time averaging. Lemma~\ref{lem:regular}
therefore gives \eqref{eq:means-general} and the equalities in
\eqref{eq:rate-general} for any periodic orbit on the carrier.
If the exponent were zero, $K_1,K_2,D$ would all vanish at one point.
Their coefficient matrix, including the constant column, is invertible,
so it cannot annihilate $(m_x,m_y,1)^{\mathsf T}$. Hence the exponent
is nonzero.

The affine functions $K_j,D-K_i$ are independent and have a common
zero, namely the actual time average. Their gradients must be
independent: dependent gradients would make the two functions
proportional once both vanish at the same point. This reasoning also
covers a zero gradient, since a constant with a zero must be identically
zero. Consequently their two equations have a unique solution.
A nonzero constant $K_i$ contradicts its cross-carrier integral, while
$K_i=0$ has already been excluded. If $m_1=m_2$, all three affine
functions vanish there, again a contradiction.
\end{proof}

The independence in Proposition~\ref{prop:means} concerns affine
functions, including their constant terms. The gradients
$\nabla K_1,\nabla K_2$ may still be parallel, and we treat that case
separately. The same argument shows that $D-K_i$ is nonconstant and
that the two lines in~\eqref{eq:means-general} are neither coincident
nor parallel. Thus their intersection is determined before any division
by a determinant. Moreover, $D$ is nonconstant: if it were constant,
Green's theorem on a periodic disk $\Omega$ would give
$\int_\Omega D\,dA=0$, hence $D=0$, contrary to affine independence.

\section{Periodic disks and an equilibrium time average}
\label{sec:periodic-disks}

We now assume the two chosen cycles are nested. The following
geometric facts are proved for all quadratic periodic orbits.

\begin{lemma}[The geometry of a quadratic periodic disk]
\label{lem:geometry}
The closed Jordan disk bounded by a quadratic periodic orbit is strictly
convex, and its time average belongs to the open disk. Every line through
an interior point meets the boundary transversely at two endpoints,
where its normal velocities have opposite signs. The open disk contains
exactly one equilibrium $e$, and $D\mathcal{X}(e)$ has two nonreal conjugate
eigenvalues. Nested periodic disks have this equilibrium in common,
and the annulus between their boundaries contains no equilibrium.
\end{lemma}
\begin{proof}
The orbit is a regular analytic Jordan curve. If its disk were not
convex, there would be two interior points with an exterior point
between them on their joining segment. A small generic perturbation
of the line preserves these three open conditions and makes its
boundary intersections transverse. Thus the line has at least four
transverse intersections with the boundary. Parameterize it by $p+tv$.
At successive intersections it alternately enters and leaves the disk;
the fixed orientation of the flow on the boundary makes
$\det(v,\mathcal{X}(p+tv))$ alternate in sign. This polynomial in $t$ has degree
at most two, whereas four intersections require at least three distinct
roots. The identically zero case cannot have a transverse intersection.
This contradiction proves convexity.

A boundary segment would make the normal velocity vanish identically
on its line. That line would be invariant, and uniqueness of solutions
at a nonstationary point would force the periodic orbit to lie in it.
Thus there is no boundary segment, which is the strict convexity meant
here; positive curvature everywhere is not required. The average of a
supporting affine functional is strictly less than its maximum, because
the orbit is not contained in the supporting line. Hence the time average
lies in the open disk. A tangent to a differentiable convex boundary
is supporting, so a line through an interior point is transverse at
both endpoints. The entry and exit orientations give the opposite signs.

If two equilibria lay inside, the normal velocity on their joining line
would have two distinct interior roots but opposite signs at the two
boundary endpoints. A polynomial of degree at most two cannot do this.
There is at least one equilibrium by the boundary index, which is one.
Thus the equilibrium is unique.
For completeness let $J(n_1,n_2)=(-n_2,n_1)$ and let $\mathcal{X}_2$ be the
quadratic homogeneous part. On the line $p+tJn$ the normal velocity is
\[
 (n\cdot \mathcal{X}_2(Jn))t^2+(n\cdot D\mathcal{X}(p)Jn)t+n\cdot \mathcal{X}(p).
\]
Opposite signs at the endpoints imply
\begin{equation}\label{eq:contact}
 [n\cdot D\mathcal{X}(p)Jn]^2
 -4[n\cdot \mathcal{X}_2(Jn)][n\cdot \mathcal{X}(p)]>0\qquad(n\ne0).
\end{equation}
If the quadratic coefficient vanishes, the linear coefficient must be
nonzero, so the same strict inequality holds. Constant polynomials,
zero polynomials and repeated roots cannot produce opposite signs.
At $p=e$ the second term vanishes. A real eigenvector $Jn$ of $D\mathcal{X}(e)$
would make the first term zero. Hence no real eigenvector exists,
and the eigenvalues are nonreal conjugates. Finally, nested disks
share the inner equilibrium by uniqueness in the outer disk.
\end{proof}

Lemma~\ref{lem:geometry} does not require the trace of
$D\mathcal{X}(e)$ to be nonzero. The next local result therefore
includes the case of purely imaginary eigenvalues.

\begin{lemma}[Local carriers and nonreal eigenvalues]
\label{lem:focus}
At a real equilibrium with nonreal conjugate eigenvalues, a real
absolutely irreducible invariant polynomial vanishing there has exactly
two smooth conjugate complex local branches with distinct tangents.
The real point is isolated on the curve, and the cofactor at that point
equals the divergence. Two distinct real absolutely irreducible
invariant curves cannot both pass through that point.
\end{lemma}
\begin{proof}
Let $F$ be the invariant polynomial and $K$ its cofactor.
Translate the point to the origin and diagonalize over $\C$ with coordinates
$(z,w)$ and eigenvalues $\lambda,\bar\lambda$, where
$\operatorname{Im}\lambda\ne0$. The real structure exchanges $z,w$
and conjugates coefficients. The lowest nonzero homogeneous term $H$
of the polynomial satisfies
\[
 (\lambda z\partial_z+\bar\lambda w\partial_w)H=K(0)H.
\]
Since $K(0)$ is real, every monomial $z^p w^q$ occurring in $H$
has $p=q$. Thus $H$ is a nonzero multiple of $(zw)^p$, with $p\ge1$,
and $K(0)=pD(0)$.

A reduced complex algebraic curve has convergent Puiseux local
branches; see \cite[Lemmas~11--12]{ChavarrigaGiacominiGrau2005Local}. Each branch is invariant,
since at its smooth points away from other branches the vector field
is tangent to the zero set. For a branch tangent to the $z$-axis write
$w=h(z)$ in Puiseux form, of order greater than one unless it is the
zero graph. With the field written as
$(A,B)=(\lambda z+O(2),\bar\lambda w+O(2))$, its equation is
$A(z,h)h'=B(z,h)$. At rational order $r>1$ the multiplier of the
new coefficient is $r\lambda-\bar\lambda$, which never vanishes.
On the fixed Puiseux denominator lattice of such a branch this
determines all integer coefficients uniquely and forces every
noninteger coefficient to vanish: at the first such order all lower
inputs have integer powers. Thus there is at most one branch in that
tangent direction, and it is smooth. The same argument applies to the other axis.
This reasoning remains valid when $\operatorname{Re}\lambda=0$.

Real conjugation forces both branches to occur. A repeated local branch
would make $F$ and both its partial derivatives vanish along that
branch, contradicting the finite singular locus of an absolutely
irreducible plane curve. Therefore $p=1$ and $K(0)=D(0)$.
The real tangent cone is definite, so the point is isolated in the
real zero set of $F$.
If two different carriers passed through it, they would share an actual
complex branch, giving infinitely many complex intersections and
contradicting B\'ezout's theorem. This is the focal local configuration
of \cite[Lemma~16]{ChavarrigaGiacominiGrau2005Local}, with the trace-zero boundary included explicitly.
\end{proof}

\begin{lemma}[A periodic average at an equilibrium]
\label{lem:meanline}
Suppose a real quadratic field with nonzero quadratic homogeneous part
has a periodic orbit whose time average is an equilibrium $e$, and
$D\mathcal{X}(e)$ has no real eigenvector. Then the field has a real invariant
affine line. In fact it is transformed by a real invertible linear
change after translation, and a positive constant time rescaling, into
\begin{equation}\label{eq:ZW}
 \xi'=-\eta+\alpha\xi+\ell\xi^2+m\xi\eta+n\eta^2,
 \qquad \eta'=\xi(1+\beta\eta),\qquad \beta\ne0.
\end{equation}
Consequently it has at most one limit cycle.
\end{lemma}
\begin{proof}
Translate $e$ to the origin and let $\langle\cdot\rangle$ denote the
time average along the given periodic orbit, in the original time.
The matrix
$C=\langle zz^{\mathsf T}\rangle$ is positive definite: if
$a^{\mathsf T}Ca=0$ for $a\ne0$, then the orbit lies in the line
$a\cdot z=0$, which is impossible. Push the field forward by
$\zeta=C^{-1/2}z$, without changing time. Now
$\langle \zeta\rangle=0$ and $\langle \zeta\zeta^{\mathsf T}\rangle=I$.
The average of each coordinate derivative is zero, as is the average
of its linear part. Thus the symmetric matrix of each quadratic
component has trace zero. Both components belong to
$\operatorname{span}_{\R}\{x^2-y^2,2xy\}$.

We choose an orthogonal rotation for which the second quadratic
component is $\mu xy$ with $\mu\ne0$. If the two quadratic components
are linearly independent, the continuous odd function
$v\mapsto\det(v,\mathcal{X}_2(v))$ on the unit circle has a zero. Take such a
$v$ as the horizontal direction. The normal quadratic component
vanishes at $(1,0)$ and is not identically zero by independence.
Its trace is zero, so both its $x^2$ and $y^2$ coefficients vanish,
as required. If the components have rank one, write $\mathcal{X}_2=Bq$, where
$B\ne0$ is constant and $q$ is a nonzero trace-free quadratic form.
Its two real zero lines are distinct. Choose one not parallel to $B$
as the horizontal direction. The same argument applies to the nonzero
normal component. No simplicity of a projective direction is assumed.

The second equation therefore has the form
\begin{equation}\label{eq:normal-mean}
 \dot y=bx+sy+\mu xy=(b+\mu y)x+sy.
\end{equation}
Here $b\ne0$, since otherwise the horizontal direction would be an
eigendirection of the linear part. Put $w=b+\mu y$. If $s\ne0$,
then on $w=0$ one has $\dot w=-sb$, a fixed nonzero constant.
A periodic function cannot have all its zeros simple with the same
crossing direction. Hence $w$ never vanishes on the orbit. It has
the sign of its average $b$. Only now do we integrate its logarithm:
\begin{equation}\label{eq:strict-average}
 0=\frac1\mu\langle \mathcal{X}\log|w|\rangle
   =\langle x\rangle+s\left\langle \frac{y}{w}\right\rangle
   =-\frac{s\mu}{b}\left\langle \frac{y^2}{w}\right\rangle.
\end{equation}
The last equality uses $y/w=y/b-\mu y^2/(bw)$ and
$\langle x\rangle=\langle y\rangle=0$. Since $w$ has fixed sign and $y$ is not identically
zero on a nonconstant periodic orbit, the last expression is nonzero.
This contradiction gives $s=0$. Equation~\eqref{eq:normal-mean} now
implies $\mathcal{X}w=\mu xw$, proving that $w=0$ is a real invariant line.

Write the first equation as $\dot x=px+hy+P_2(x,y)$. The nonreal
spectrum implies $-hb=\det D\mathcal{X}(0)>0$. Set
$\omega=\sqrt{-hb}$, $\tau=\omega t$, $\xi=bx/\omega$, and
$\eta=y$. In the new time $\tau$ the equations are precisely
\eqref{eq:ZW}, with $\alpha=p/\omega$ and $\beta=\mu/b\ne0$.
This real invertible coordinate change and positive constant time
rescaling preserve periodic orbits and isolation. Theorem~1 of
Zhou and Wang \cite[pp.~370--371, system~(1.1)]{ZhouWang2003Line} gives at most
one limit cycle for this normal form.
\end{proof}

The covariance in this proof is that of the given periodic orbit.
The conclusion also uses its zero first moments and the periodic
logarithmic integral in~\eqref{eq:strict-average}; a positive definite
solution of the quadratic moment equations alone is insufficient.
The orbit in Lemma~\ref{lem:meanline} need not be isolated.

\section{The common equilibrium and parallel cofactors}
\label{sec:parallel-cofactors}

We return to the chosen nested pair, and let $e$ be the common
interior equilibrium supplied by Lemma~\ref{lem:geometry}.

\begin{proposition}[The common equilibrium lies outside the carriers]
\label{prop:avoid-focus}
The two carriers avoid $e$. In particular,
\begin{equation}\label{eq:common-focus}
 F_1(e)F_2(e)\ne0,\qquad K_1(e)=K_2(e)=0,\qquad D(e)\ne0.
\end{equation}
\end{proposition}
\begin{proof}
By Lemma~\ref{lem:focus}, at most one carrier passes through $e$.
If $F_i(e)=0$, that lemma gives $K_i(e)=D(e)$, while invariance of
the nonvanishing other carrier gives $K_j(e)=0$. The unique solution
in Proposition~\ref{prop:means} then gives $e=m_i$.
Lemma~\ref{lem:meanline} applies to $\Gamma_i$ and implies that the
system has at most one limit cycle, a contradiction. Thus both
carriers avoid $e$. Evaluating their invariant equations at $e$ gives
$K_1(e)=K_2(e)=0$. If $D(e)=0$ as well, the three independent affine
functions would have a common zero. Therefore $D(e)\ne0$, and $e$
is a strong focus.
\end{proof}

\begin{proposition}[Exclusion of parallel cofactor gradients]
\label{prop:parallel}
The gradients of the two cofactors cannot be parallel in a nested pair.
\end{proposition}
\begin{proof}
Both cofactors are nonconstant by Proposition~\ref{prop:means}.
Parallel gradients give $K_2=cK_1+d$ with $c\ne0$. Independence
of the two affine functions gives $d\ne0$.
If $c<0$, there is a direct exclusion, valid even without nesting.
The function $H=|F_1|^{-c}|F_2|$ is continuous on the plane, zero
on the two carriers, and smooth off them, where $\mathcal{X}H=dH$.
On the closed disk bounded by $\Gamma_1$, it is zero on the boundary
and positive at some interior point, since a nonzero polynomial has
zero set with empty interior. Its positive maximum is attained at an
interior point off both carriers, where $\mathcal{X}H=0$, contradicting $dH\ne0$.
Isolated real zeros and other interior branches do not affect this
maximum argument.

If $c>0$, Proposition~\ref{prop:avoid-focus} gives
$0=K_2(e)=cK_1(e)+d=d$, again a contradiction.
The role of Lemma~\ref{lem:meanline} in this case can also be seen
directly. The relation $d\ne0$ prevents both cofactors from vanishing
at $e$, so invariance and Lemma~\ref{lem:focus} require exactly one
carrier to pass through $e$. Its time average is then $e$.
Indeed,
$1,K_1,D$ are independent, so $(u,v)=(K_1,D)$ are valid affine
coordinates. The mean equations give
\[
 (u,v)(m_1)=(-d/c,-d/c),\qquad (u,v)(m_2)=(0,d).
\]
The carrier through $e$ has its corresponding mean there, and
Lemma~\ref{lem:meanline} excludes the pair. This argument does not
require the two chosen cycles to be adjacent.
\end{proof}

The affine coordinates used above preserve time. The only time
rescaling is the positive constant one at the end of
Lemma~\ref{lem:meanline}.

\section{Independent cofactor gradients and opposite stability}
\label{ind:sec:means}

It remains to consider
$\det(\nabla K_1,\nabla K_2)\ne0$.
The preceding reduction gives a common interior equilibrium $e$ with
\begin{equation}\label{ind:focus}
 F_1(e)F_2(e)\ne0,\qquad K_1(e)=K_2(e)=0,
 \qquad \delta:=D(e)\ne0.
\end{equation}
In particular, $e$ is a strong focus. We use the real affine coordinates
$(u,v)=(K_1,K_2)$, so that $e=(0,0)$. If this coordinate map is
$\Phi(z)=Mz+h$, the transformed field is $M\mathcal{X}\circ\Phi^{-1}$ and
its derivative is $M\,D\mathcal{X}(\Phi^{-1})M^{-1}$. Thus its divergence is
$D\circ\Phi^{-1}$. No change of time is involved: periods and
characteristic exponents are unchanged, while time averages transform
by $\Phi$. We use
the same letters for the transformed quantities. There are unique
$a,b\in\R$ such that
\[
 D=\delta+au+bv.
\]
Proposition~\ref{prop:means} gives
\begin{equation}\label{ind:means}
 m_1=(s_1,0),\qquad m_2=(0,s_2),\qquad
 (1-a)s_1=\delta,\qquad (1-b)s_2=\delta.
\end{equation}
Since $\delta\ne0$, equations~\eqref{ind:means} imply that $1-a$,
$1-b$, $s_1$ and $s_2$ are nonzero.
Every periodic orbit on $F_i=0$ has the same mean $m_i$ and
the same characteristic exponent per unit time $s_i$; the cofactor
integral formulas apply even if that orbit was not initially assumed
isolated. In particular, all such orbits are hyperbolic.

Let $\Omega$ be the open disk bounded by the outer chosen cycle.
The real affine line
\begin{equation}\label{ind:S}
 S:=D-K_1-K_2=0
 \quad\Longleftrightarrow\quad
 \frac{u}{s_1}+\frac{v}{s_2}=1
\end{equation}
contains $m_1$ and $m_2$. By Lemma~\ref{lem:geometry}, the entire
closed segment joining these two distinct means lies in $\Omega$.
Notice also that $S(e)=\delta\ne0$.

\begin{lemma}[Opposite stabilities]
\label{ind:sign}
There is no periodic orbit outside the union of the two carriers,
and $s_1s_2<0$.
\end{lemma}
\begin{proof}
A periodic orbit meeting an invariant carrier is contained in it,
by uniqueness of solutions. If a periodic orbit avoids both carriers,
integration of $\mathcal{X}\log|F_i|=K_i$ along it shows that the averages of
both $u$ and $v$ vanish. Its time average is therefore $e$.
Since the eigenvalues at $e$ are nonreal, Lemma~\ref{lem:meanline}
implies that the system has at most one limit cycle, contradicting
the chosen pair. This argument does not require the extra
periodic orbit to be isolated.

Suppose now that $s_1s_2>0$. Every periodic orbit lies on one of the
two carriers and is hyperbolic with the same stability. Choose the
direction of time in which all these orbits repel. The closed annulus
between the chosen cycles is compact, invariant and contains no
equilibrium, by Lemma~\ref{lem:geometry}. In its interior choose a
point where $F_1F_2\ne0$; such points exist because a nonzero
polynomial has zero set with empty interior. Its orbit is not
periodic, by the first assertion. Its positive semiorbit remains
in the compact annulus, so the Poincar\'e--Bendixson theorem implies
that its $\omega$-limit set is a periodic orbit. A nonperiodic positive
semiorbit cannot tend to a hyperbolic repelling cycle, whether that
cycle lies in the interior or on the boundary of the annulus. This
follows directly from the expanding return map on either side of
such a cycle. This is a contradiction. Since neither $s_i$ is zero,
their product is negative.
\end{proof}

The argument applies even when other cycles lie between the two
chosen cycles. The use of cofactor averages to determine
stability has precedents in Gasull and Giacomini
\cite[Section~5, Corollary~5.3]{GG2023}; here the exclusion of periodic
orbits outside the carriers additionally uses the periodic-average
argument of Lemma~\ref{lem:meanline}.

\section{Cofactor zero lines and local invariant branches}
\label{ind:sec:local}

We first isolate the two consequences of the projective two-point
criterion that we need. In its published form
\cite[Theorem~2.4]{ChavarrigaGarcia2006Cubic}, the criterion concerns a real quadratic
field with coprime components and a real invariant
algebraic curve of degree at least two, with cofactor $K$. Write
\[
 P^{h}(U,V,Z)=Z^2P(U/Z,V/Z),\qquad
 Q^{h}(U,V,Z)=Z^2Q(U/Z,V/Z),
\]
and $K^{h}(U,V,Z)=ZK(U/Z,V/Z)$. If two distinct points of
$\mathbb{CP}^{2}$ satisfy $P^{h}=Q^{h}=K^{h}=0$, all limit
cycles are contained in that curve. In the present setting $e$
already supplies one of the two points for either cofactor.
For completeness, the following elementary proof covers both
uses of this criterion without invoking its general statement.

\begin{lemma}[Zeros of a cofactor]
\label{ind:line}
Let $F$ be either carrier and $K$ its nonconstant affine cofactor.
Under the preceding assumptions, the line $K=0$
contains no finite complex equilibrium other than $e$. Moreover,
the two components of $\mathcal{X}_2$ cannot vanish simultaneously in the
direction of this line.
\end{lemma}
\begin{proof}
Use real affine coordinates with $x=K$ and with $e$ as the origin
on $x=0$. If a second finite complex equilibrium lies on this
line, the restrictions of both components of $\mathcal{X}$ to the line
are real polynomials of degree at most two with the same two
distinct complex roots. These two restrictions span a real vector
space of dimension at most one. Indeed, if the second root is real,
both restrictions are multiples of the same quadratic; if it is
nonreal, reality supplies its conjugate as a third distinct root,
so both restrictions vanish identically. If instead $\mathcal{X}_2$ vanishes
in the direction of the line, the two restrictions have degree at
most one and share the root $e$, giving the same dimension bound.

In either case there is a nonzero real linear form $\ell$ for
which $\mathcal{X}\ell$ vanishes identically on $x=0$. If its linear part
is proportional to that of $x$, then $x=0$ is invariant. The cycle
on the other carrier cannot be contained in or cross this line,
so $x$ has a fixed nonzero sign along that cycle. This contradicts
the cross-carrier identity $\int x\,dt=0$. If both component
restrictions vanish identically, this same argument applies.

Otherwise choose $y=\ell$, completing $x$ to affine coordinates.
Divisibility by $x$ gives
\begin{equation}\label{ind:line-form}
 \mathcal{X}y=x(Ax+By+C),\qquad A,B,C\in\R.
\end{equation}
The polynomial $F$ is nonzero on a collar
of the cycle on the other carrier. For $B\ne0$, define there
\[
 G=(y+C/B)|F|^{-B}.
\]
Since $\mathcal{X}F=xF$, direct differentiation yields
\begin{equation}\label{ind:line-integral}
 \mathcal{X}G=A x^2|F|^{-B}.
\end{equation}
For $B=0$, instead put $G=y-C\log|F|$; then $\mathcal{X}G=A x^2$.
Integrating either identity along the cycle gives $A=0$, because
the positive weight is well-defined and the orbit is not contained
in $x=0$.

Consequently $G$ is a real analytic first integral on the collar.
It cannot be constant on any nonempty open subset. For $B\ne0$,
constancy, after restriction to an open subset with $y+C/B\ne0$
and differentiation with respect to $x$, would imply $F_x=0$
there, hence identically. For $B=0$ and $C\ne0$ the same follows
by differentiation. A polynomial depending only on $y$ cannot
contain a nonconstant periodic orbit, since a continuous orbit
in its zero set would have constant $y$. The remaining case
$B=C=0$ gives $G=y$, which is manifestly nonconstant.

On a transverse analytic section the restriction of $G$ is also
nonconstant: otherwise flow invariance and a flow box would make
it constant on an open set. Its analytic restriction is strictly
monotone on a sufficiently short punctured side of the cycle.
The orientation-preserving return map preserves that restriction
and is therefore the identity on that side. Analyticity then
excludes isolation of the cycle, the required contradiction.
\end{proof}

It follows that every finite complex equilibrium $q\ne e$ satisfies
\begin{equation}\label{ind:both}
 K_1(q)K_2(q)\ne0,
 \qquad F_1(q)=F_2(q)=0.
\end{equation}
The second assertion follows by evaluating $\mathcal{X}F_i=K_iF_i$ at $q$.

The next lemma gives the local branch properties needed below;
compare Chavarriga, Giacomini and Grau
\cite[Theorem~7(a), Lemmas~11--12]{ChavarrigaGiacominiGrau2005Local}.
We give the argument for algebraic branches in the stated spectral ranges.

\begin{lemma}[Local branches and cofactor values]
\label{ind:branches}
Let $\mathcal{X}$ be an analytic vector field vanishing at a complex
point $q$, with linearization eigenvalues $\lambda,\nu$. If a nonzero
polynomial $F$ satisfies $F(q)=0$ and $\mathcal{X}F=KF$, then
\begin{equation}\label{ind:weight-spectrum}
 K(q)=p\lambda+r\nu
 \quad\text{for some }p,r\in\Z_{\ge0},\qquad
 p+r=\operatorname{ord}_q F.
\end{equation}
Suppose, in addition, either that $\lambda\nu\ne0$ and
$\lambda/\nu\notin\mathbb Q_{>0}$, or that exactly one
eigenvalue is zero. Then at most one actual algebraic invariant
branch has either eigendirection as tangent, and any such branch
is smooth. If two distinct absolutely irreducible invariant
curves pass through $q$, each has exactly one local branch there;
their cofactor values are the two respective normal eigenvalues
and sum to the trace of the linear part.
\end{lemma}
\begin{proof}
Translate $q$ to zero, put $L=D\mathcal{X}(q)$, and let $H$ be the first
nonzero homogeneous part of $F$, of degree $d$.
The lowest-degree part of the invariance equation is
\[
 (Lz)\cdot\nabla H=K(q)H.
\]
After triangularizing $L$, the operator on homogeneous polynomials
of degree $d$ is triangular in a monomial basis, with diagonal
entries $p\lambda+r\nu$, $p+r=d$. This proves
\eqref{ind:weight-spectrum}, including repeated eigenvalues,
nondiagonalizable matrices and the double-zero spectrum.

Under either of the additional spectral assumptions the two
eigenvalues are distinct, so we may use complex eigenvector
coordinates. An irreducible algebraic branch admits a convergent
Puiseux parametrization. It is invariant: away from its intersections
with the other local branches, the identity $\mathcal{X}F=KF$ makes the
field tangent to its smooth part. Its tangent is an eigendirection.
Indeed, writing a parametrization as $z(t)=t^h v+\cdots$ and
comparing the leading terms of tangency gives
$\det(v,Lv)=0$, including the possibility $Lv=0$.

For a branch tangent to the horizontal axis write a Puiseux graph
$y=f(x)$, with order greater than one (the zero graph is allowed).
Its equation is
\begin{equation}\label{ind:puiseux}
 (\lambda x+O(2))f'(x)=\nu f(x)+O(2),
\end{equation}
where the higher-order terms are evaluated on the graph. At any
rational exponent $r>1$, the coefficient yet to be determined
has multiplier $r\lambda-\nu$. This is nonzero when both
eigenvalues are nonzero and their ratio is not a positive rational
number; it is also nonzero if exactly one eigenvalue is zero.
Recursion on a fixed Puiseux denominator lattice therefore gives
at most one series. The unique integer-power formal solution
lies on every such lattice, so all noninteger coefficients of an
actual Puiseux solution vanish. Hence an actual branch is a smooth
graph and is unique for this tangent. Interchanging the axes
proves the same assertion for the other tangent. Only branches supplied
by an algebraic curve are used, so convergence of a general formal
weak separatrix is not required.

A reduced absolutely irreducible algebraic curve cannot contain
a repeated local branch: otherwise it and its first partial
derivatives would vanish on that branch, producing infinitely
many singular points. Nor can two distinct absolutely irreducible
curves share an actual complex branch, since their finite
intersection would then contain infinitely many points. Thus two
such curves through $q$ must occupy the two different smooth
branches, one each. Their first terms are nonzero multiples of
the respective normal coordinates. The lowest-degree equation
then gives the normal eigenvalue as the cofactor value, proving
the last assertion.
\end{proof}

Nodes with positive rational eigenvalue ratio fall outside this branch
classification; below we use only their positive index.

\section{Finite equilibria}
\label{ind:sec:finite}

\begin{proposition}[Finite equilibria and their indices]
\label{ind:finite}
Every finite complex equilibrium is a real simple equilibrium,
and its real local index is $+1$.
\end{proposition}
\begin{proof}
We first exclude nonreal equilibria. Let $q\notin\R^2$, with
linearization eigenvalues $\lambda,\nu$. If these two complex
numbers are linearly dependent over $\R$, their real span has
dimension at most one. By \eqref{ind:weight-spectrum}, this span
contains $K_1(q)$ and $K_2(q)$, as well as
$D(q)=\lambda+\nu$. It therefore contains
\[
 D(q)-aK_1(q)-bK_2(q)=\delta\in\R\setminus\{0\}.
\]
The span must be exactly the real axis. Both $K_i(q)$ are real;
since $(K_1,K_2)$ is an invertible real affine coordinate map,
$q$ is real, a contradiction. This argument includes zero
eigenvalues, repeated eigenvalues and positive rational ratios.

If instead $\lambda,\nu$ are linearly independent over $\R$,
their ratio is nonreal. By \eqref{ind:both}, both carriers pass
through $q$, and Lemma~\ref{ind:branches} gives
$D(q)=K_1(q)+K_2(q)$. Hence $q$, and by reality also $\bar q$,
lie on the complexification of the real line $S=0$ in
\eqref{ind:S}. This line cuts the interior of the outer periodic
disk. Parametrize it with a real affine parameter. The normal
component of $\mathcal{X}$ restricted to it is a real polynomial of degree
at most two. At the two boundary endpoints of the chord it has
opposite nonzero signs, by Lemma~\ref{lem:geometry}, and so it
has a real root between them. The distinct nonreal points $q$
and $\bar q$ provide two additional nonreal roots. This is
impossible for a polynomial of degree at most two; it cannot
vanish identically because the endpoint values are nonzero.

Now let $q$ be real and different from $e$. If both eigenvalues
are zero, \eqref{ind:weight-spectrum} gives $K_1(q)=K_2(q)=0$,
contrary to \eqref{ind:both}. If exactly one eigenvalue is zero,
Lemma~\ref{ind:branches} makes the two cofactor values the two
normal eigenvalues, one of which is zero, again a contradiction.
Thus $D\mathcal{X}(q)$ is invertible. The focus $e$ is also nondegenerate.
All finite equilibria are consequently simple zeros of the
polynomial map $(P,Q)$.

If $q\ne e$ were a saddle, its eigenvalues would be real,
nonzero and of opposite signs. Lemma~\ref{ind:branches} would
give $S(q)=0$ and $K_1(q)K_2(q)<0$. There is a unique real
$t$ such that
\[
 q=((1-t)s_1,ts_2),\qquad
 K_1(q)K_2(q)=t(1-t)s_1s_2.
\]
Because $s_1s_2<0$, the strict negativity of this product is
equivalent to $0<t<1$. Thus $q$ lies in the open segment between
the two means, and hence in $\Omega$. But the only equilibrium
in $\Omega$ is $e$, and $S(e)=\delta\ne0$. This contradiction
excludes saddles. The determinant of the linearization at every
real simple equilibrium is therefore positive, giving local index
$+1$. In particular, real nodes with positive rational eigenvalue
ratio are included without imposing any bound on their local
invariant branches.
\end{proof}

\section{Common directions at infinity and the final count}
\label{ind:sec:infinity}

We retain $e$ as the origin, but may make a further real linear
coordinate change. Write $P_j,Q_j$ for the homogeneous parts
of the two components and $k_i$ for the degree-one homogeneous
part of $K_i$. In these coordinates $K_i=k_i$, since $K_i(e)=0$.

\begin{lemma}[Common directions at infinity]
\label{ind:infinity}
The homogeneous quadratic polynomials $P_2,Q_2$ have no common
zero in $\mathbb{CP}^{1}$.
\end{lemma}
\begin{proof}
Suppose there is a common direction $[1:t_0]$, with $t_0$
possibly complex. Interchanging the coordinates if necessary
covers a vertical direction. Put
\[
 p(t)=P_2(1,t),\qquad H(t)=Q_2(1,t)-tp(t).
\]
For either carrier, denote its defining polynomial and cofactor by
$F$ and $K$, and write $k$ for the linear homogeneous part of $K$.
Set $n=\deg F$ and
$F_n(x,y)=x^n f(y/x)$. Comparison of terms of degree $n+1$
in $\mathcal{X}F=KF$ gives the polynomial identity
\begin{equation}\label{ind:top}
 H(t)f'(t)+np(t)f(t)=k(1,t)f(t).
\end{equation}
At $t_0$ we have $p(t_0)=H(t_0)=0$.
Writing $f(t)=(t-t_0)^r g(t)$ with $g(t_0)\ne0$ and
$r\ge0$, division by $(t-t_0)^r$ in this identity followed
by evaluation at $t_0$ yields
\begin{equation}\label{ind:root}
 k(1,t_0)=rH'(t_0).
\end{equation}
For $r=0$, this follows directly from \eqref{ind:top}.
In particular, if $H'(t_0)=0$, both $k_i$ vanish in this
direction. This includes $H\equiv0$, so the radial quadratic
case requires no branch argument along a nonisolated singular
set at infinity.

Suppose that $H'(t_0)\ne0$ and that both $k_i(1,t_0)$ are
nonzero. Equation~\eqref{ind:root} implies $r_i>0$, so both
projective carrier closures pass through $[1:t_0:0]$.
In the projective chart
\[
 t=y/x,\qquad w=1/x,
\]
set
\[
 A(t,w)=w^2P(1/w,t/w),\qquad
 B(t,w)=w^2Q(1/w,t/w).
\]
Multiplying the coordinate expression of the original field by
$w$ gives the polynomial local field
\begin{equation}\label{ind:chart-field}
 \mathcal{Y}=(B-tA)\partial_t-wA\partial_w.
\end{equation}
The line $w=0$ is invariant. At $(t_0,0)$ the linearization
is upper triangular, with diagonal entries $H'(t_0)$ and
$-p(t_0)=0$. The invariant line $w=0$ is tangent to the eigendirection
corresponding to the nonzero eigenvalue.

The local equation of the projective closure of a carrier is
$\widehat F=w^nF(1/w,t/w)$, and direct differentiation gives
\begin{equation}\label{ind:chart-cofactor}
 \mathcal{Y}\widehat F=
 \bigl(wK(1/w,t/w)-nA(t,w)\bigr)\widehat F.
\end{equation}
Thus both closures supply actual invariant algebraic branches
of $\mathcal{Y}$. By Lemma~\ref{ind:branches}, apart from $w=0$ there
is at most one further actual branch at this semihyperbolic
point. Neither projective carrier contains the line at infinity,
and the two distinct projective irreducible curves cannot share
an actual branch, by B\'ezout's finite-intersection property.
They therefore cannot both pass through this point. We conclude
that, also when $H'(t_0)\ne0$, at least one $k_i(1,t_0)$
must vanish.

Multiplication by $w$ is confined to the local projective analysis;
all periodic averages above use the original time. The degree
correction in~\eqref{ind:chart-cofactor} vanishes at $(t_0,0)$ because
$A(t_0,0)=p(t_0)=0$.

Finally, a nonzero real linear form $k_i$ has a real projective
zero direction. Hence the common direction just obtained is real
and is the direction of the line $K_i=0$ through $e$.
Both $P_2$ and $Q_2$ vanish there, contradicting
Lemma~\ref{ind:line}. In terms of the published projective
criterion, the two distinct points are $e$ and this point at
infinity, and both highest homogeneous components vanish
at the latter. A characteristic direction
$\mathcal{X}_2(v)\parallel v$ would not suffice.
\end{proof}

The use of the invariant line at infinity to occupy one branch
at a semihyperbolic point, together with the cofactor correction
in \eqref{ind:chart-cofactor}, has a direct precedent in
Sorolla \cite[Chapter~5, Lemma~5.6]{Sorolla2004Thesis}. There it appears
under a non-nested configuration. In the present argument the
common finite focus and Lemma~\ref{ind:line} supply the required
finite point for each cofactor.

\begin{proposition}[Exclusion of independent cofactor gradients]
\label{ind:impossible}
The cofactor gradients of the chosen nested pair cannot
be linearly independent.
\end{proposition}
\begin{proof}
By Lemma~\ref{ind:infinity}, $P_2$ and $Q_2$ have no common
complex projective zero. In particular, neither is identically
zero: if one were zero, the other, a nonzero binary homogeneous
quadratic, would have a zero in $\mathbb{CP}^{1}$. Thus both
components have degree two. Their homogenizations have no
common component. Indeed, an affine common component is excluded
by coprimality of $P,Q$, and the line at infinity is not
a common component because both leading terms are nonzero.
B\'ezout's theorem therefore gives total complex intersection
number four. There are no intersections at infinity, so all
four lie in the affine plane. By Proposition~\ref{ind:finite}, every
affine intersection is real and simple. The nonzero Jacobian determinant
implies that the two component curves meet transversely, with
intersection multiplicity one. There are exactly four distinct
finite equilibria, each with real index $+1$.

The index of $\mathcal{X}$ on a sufficiently large positively oriented
circle is consequently $4$. On the other hand, $\mathcal{X}_2$ is nowhere
zero on the real unit circle, so its norm has a positive lower
bound there. On a sufficiently large circle the straight-line
homotopy which removes the lower-degree terms of $\mathcal{X}$ never
vanishes. Its index equals the winding number of $\mathcal{X}_2$ on the
unit circle. For fixed real vectors $U,V,W$, this homogeneous
quadratic field has the form
\[
 \mathcal{X}_2(\cos\theta,\sin\theta)
   =U+V\cos(2\theta)+W\sin(2\theta).
\]
If $V,W$ are linearly independent, the image is an ellipse,
traversed twice, which avoids the origin. Its winding number
is $0$, $2$ or $-2$. If $V,W$ are dependent, the image is a
line segment or a point avoiding the origin, and its winding
number is zero. Neither possibility gives $4$, a contradiction.
\end{proof}

\section{Completion of the proof}
\label{sec:completion}

\begin{proof}[Proof of Theorem~\ref{thm:main}]
Suppose that system~\eqref{eq:system} has two nested limit cycles
contained in distinct minimal carriers. The reductions in
Section~\ref{sec:preliminaries} exclude affine fields and fields with
a nonconstant common factor.
Lemmas~\ref{lem:regular} and~\ref{lem:weights}, and
Proposition~\ref{prop:means}, exclude zero or constant cofactors,
dependent cofactors and divergence in their span, while also deriving
hyperbolicity rather than assuming it. Both remaining cofactor
gradients are nonzero. Their parallel case is excluded by
Proposition~\ref{prop:parallel}, and their independent case by
Proposition~\ref{ind:impossible}. These cases exhaust the possibilities.
\end{proof}

The argument places no restriction on carrier degree and retains all
local algebraic branches. Its use of a compact invariant annulus does
not require the two chosen cycles to be adjacent. The only conclusion
about a single irreducible carrier in Corollary~\ref{cor:carriers} is
that no other minimal carrier can contain a limit cycle; the number
of limit cycles on that one carrier is not determined here.

\raggedbottom
\section*{Acknowledgements}

The first author is partially supported by the China Scholarship Council
(No. 202506230080) and the National Natural Science Foundation of China
(No. 12271355). The second author is partially supported by the Agencia
Estatal de Investigaci\'{o}n of Spain grant \mbox{PID2025-168432NB-I00}.

\section*{Data availability statement}

This paper has no data.

\section*{Author contribution statement}

Both authors contributed equally to this paper.

\end{document}